\documentclass[14t]{article}
\usepackage{times}
\newtheorem{definition}{Definition}

\newtheorem{theorem}[definition]{Theorem}
\newtheorem{proposition}[definition]{Proposition}
\newtheorem{corollary}[definition]{Corollary}
\begin{document}
\global\def\refname{{\normalsize \it References:}}
\baselineskip 12.5pt
%
%
% TITLE, AUTHOR, ABSTRACT, KEYWORDS
%
\title{\LARGE \bf On Rigid,   Hard and Soft   Problems and Results in Arithmetic Geometry}

\date{}

\author{\hspace*{-10pt}
\begin{minipage}[t]{2.7in} \normalsize \baselineskip 12.5pt
\centerline{NIKOLAJ GLAZUNOV}
\centerline{National Aviation University}
\centerline{Department of Electronics}
\centerline{1 Cosmonaut Komarov Avenue, 03680 Kiev}
\centerline{UKRAINE}
\centerline{glanm@yahoo.com}
\end{minipage} \kern 0in
%\begin{minipage}[t]{2.7in} \normalsize \baselineskip 12.5pt
%\centerline{SECOND AUTHOR}
%\centerline{Name of the University}
%\centerline{Institute of Mechanical Engineering}
%\centerline{47 West Lincoln Avenue, 87 115 City}
%\centerline{COUNTRY}
%\centerline{second.author@math.univ.ab}
%\end{minipage}
%
% If you are three authors then you can use three mini--pages
% instead of two. Their horizontal size must be less than 2.7in
% indicated above. It can be e.g. 2.3in. However, you must pay
% attention that you do not exceed the total width of the text.
%
\\ \\ \hspace*{-10pt}
\begin{minipage}[b]{6.9in} \normalsize
\baselineskip 12.5pt {\it Abstract:}
% The text of the abstract follows.
Rigid, hard and soft problems and results in arithmetic geometry are presented. "Soft" and "hard" in our paper are limited to the framework
 of solutions of quadratic forms over rings of integers of local and global fields, the Hardy-Littlewood-Kloosterman method. Next we consider the notion
 of rigidity. In the framework we give review of some novel results in the aria.
\\ [4mm] {\it Key--Words:}
% The key-words follow.
Diophantine equation, rigidity, Dirichlet character, ergodic method, hermition, Hardy-Littlewood-Kloosterman method, convolution of measures,
group action, uniform rigidity, superrigidity
\end{minipage}
\vspace{-10pt}}

\maketitle

\thispagestyle{empty} \pagestyle{empty}
% numbers of pages are supplemented by the editor
%
% THE BEGINNING OF THE TEXT
%
\section{Introduction}
\label{S1} \vspace{-4pt}

We review some novel results and methods on rigidity. These include (but not exhaust) 
methods and results by H. Furstenberg, G.A. Margulis,  G. D.  Mostow, J. Bourgain, A. Furman, A. Lindenstrauss, S. Mozes, J. James, T. Koberda, K. Lindsey, C. Silva, P. Speh, A. Ioana, K. Kedlaya, J.  Tuitman \cite{KT,BFLM,Io,JKLSS}, and others.
M. Gromov\cite{Gr} in his talk at the International Congress of Mathematicians in Berkeley have presented problems and results of soft and hard symplectic geometry. In this connection  
we will present some soft and hard problems and results in arithmetic geometry. "Soft" and "hard" in our talk are limited to the framework of solutions of algebraic equations over rings of integers of local and global fields and the elements of Hardy-Littelwood-Kloosterman methods.

It is well known that the projective space is rigid~\cite{Hartshorne}. The set of integer solutions of a Diophantine equation is a hard or a rigid object.
Lattices in spaces are rigid objects. 

\section{Sums of squares}
\label{S2} \vspace{-4pt}
Let $\Lambda$ be a lattice  in $n$-dimensional real euclidean space that is defined by congruences. Davenport, Mordell, Cassels and others used the lattices and Minkowski`s convex body theorem for proving results about existence of nontrivial solutions of some Diophantine equations. We will give examples below.
 Recall the case of positive quadratic forms.
Let $\tau$ be a complex number, $Im \; \tau > 0, q = \exp{\pi i \tau}, \theta_{3}(\tau) = \sum_{m = -\infty}^{\infty} q^{m^2}$
the Jacobi function.
Denote by ${\bf Z}^n$ the $d-$~dimensional integer lattice.
Let $r_{n}(m)$ be the number of ways of writing $m$ as a sum $f(x_1 \cdots x_n ) = f$ of $n$ squares.
Put $\Theta_{{\bf Z}^n} = \theta_{3}(\tau)^n$.

\subsection{Sums of two squares}
\vspace{-4pt}
Let $p \equiv 1 \pmod {4}$. In the case there is the integer $l$ such that $l^2 + 1 \equiv 0 \pmod {p}$. The lattice 
$\Lambda$ of pairs $(a, b)$ of integer numbers is defined by congruences $a \equiv lb \pmod {p}$ and has determinant
 $d(\Lambda) \le p$.
 From this and Minkowski`s convex body theorem follow that every prime  $p \equiv 1 \pmod {4}$ is the sum of two squares.

Let $\chi$ be the nontrivial Dirichlet character $\bmod \; 4$, integer $m > 0$. There is the well known 
\begin{proposition} \label{T0}
The number of integer solutions  of the equation $x_1^2 + x_2^2 = m $ is equal $4 \sum_{d|m} \chi(d)  $.

In the framework of the function  $\Theta_{{\bf Z}^n}$   we have  
 $\Theta_{{\bf Z}^2} =  \sum_{m = 0}^{\infty} r_{2}(m) q^{m} $.
\end{proposition}

\subsection{Sums of three squares}
\vspace{-4pt}
In the case and in the case $n = 4$ it is possible to use quaternions (hermitions) but for simplicity we will
 formulate the well known result by $\Theta_{{\bf Z}^3}$ and $r_{3}(m)$.
\begin{proposition} \label{T01}
$\Theta_{{\bf Z}^3} =  \sum_{m = 0}^{\infty} r_{3}(m) q^{m} $.
\end{proposition}

\subsection{Sums of four squares}
\vspace{-4pt}
The quadratic form $x_1^2 + x_2^2 + x_3^2  + x_4^2   $ represents all positive numbers (Lagrange). The number of solutions 
of the equation 
$x_1^2 + x_2^2 + x_3^2  + x_4^2   = m $, where $m$ is a positive integer is given by Jacobi.
\begin{proposition} \label{T1}
The number of integer solutions  of the equation $x_1^2 + x_2^2 + x_3^2  + x_4^2   = m $ is equal $8 \sum_{d|m} d  $ if $m = 2k + 1$,
and is equal $24 \sum_{d|m} d  $ if $m = 2k .$

In the framework of the function  $\Theta_{{\bf Z}^n}$   we have
 $\Theta_{{\bf Z}^4} =  \sum_{m = 0}^{\infty} r_{4}(m) q^{m} $.
\end{proposition}

\subsection{Sums of  squares greater than four}
\vspace{-4pt}
Recall elements of Hardy-Littelwood-Kloosterman method in the case.
This is valid also in the previous case $n = 4$. Consider a function
 of complex variable $u, |u| < 1$
$$ \vartheta (f,u) =   \sum_{x_1 \cdots x_n = -\infty}^{\infty} u^{f(x_1 \cdots x_n )}   $$
Then the number $r_{n}(m)$ of ways of writing $m$ as a sum of $n$ squares
by Cauchy's integral formula is given as
$$r_{n}(m) = \frac{1}{2\pi i}\oint_{\Gamma} \vartheta (f,u) u^{-m - 1} du   $$
where $\Gamma$ is the circle $|u| =  \exp(-\frac{1}{m})$.
We omit here the very important step of the dividing $\Gamma$  into Farey-arcs. 

\section{Elements of history of rigidity}
\label{S3} \vspace{-4pt}

The history of rigidity is reflected in papers by A. Selberg, E. Calabi, E. Vesentini, A. Weil, H. Furstenberg, G. Mostow, G. A. Margulis and their colleagues \cite{Se,CV,We,Fu,Mo,Ma}. There are interesting surveys by  D. Fisher \cite{Fi} and R. Spatzier  \cite{Sp}.
Let $G$ be a finitely generated group, $D$ a topological group, and $h: G \to D$ a homomorphism. Follow to \cite{Fi}  recall 

\begin{definition} \label{D1}
Given a homomorphism $h: G \to D$, it is said that $h$ is {\it locally rigid} if any other homomorphism $h^{`}$
 which is close to $h$ is conjugate to $h$ by a small element of $D$.
\end{definition}

Recall follow to  \cite{CV,We} in framework of  \cite{Sp} the Local Rigidity Theorem. 

\begin{theorem} \label{T12}
Cocompact discrete subgroups $H$ in semisimple Lie groups without compact nor $SL(2,{\bf R})$ nor 
$SL(2,{\bf C})$ local factors is deformation rigid.
\end{theorem}

The notion of uniform rigidity was introduced as a topological version of rigidity by S. Glasner and D. Maon \cite{GM}.

\section{ Uniformly rigid and  measurable weak mixing}
\label{S4} \vspace{-4pt}

Authors of the paper \cite{JKLSS}  investigate properties of uniformly rigid transformations and analyze the compatibility of uniform rigidity
 and measurable weak mixing along with some of their asymptotic convergence properties.
\par
This interesting survey includes some resent results on genericity of rigid and multiply recurrent infinite measure
 preserving and nonsingular transformations by O. Ageev and C. Silva \cite{AS} and on measurable sensitivity by J. James, T. Koberda, K. Lendsey, C. Silva,
  P.~Speh  \cite{JKLSS1}. 
All spaces of the paper \cite{JKLSS} are considered simultaneously as topological spaces and as measure spaces.
 Presented results concern either the measurable dynamics on the spaces or the interplay between the measurable
 and topological dynamics. 
After some introductory section, second section of the paper\cite{JKLSS}  considers functional analytic properties
 of uniform rigidity that is similar to the properties of rigidity. Authors formulate and prove
\begin{theorem} \label{T13}
Theorem 1. Every totally ergodic finite measure-preserving transformation on a Lebesgue space has a representation
 that is not uniformly rigid, except in the case where the space consists of a single atom.
\end{theorem}
\par
The proof of the theorem connects with results of authors of the paper \cite{JKLSS} that uniform rigidity and weak
 mixing are mutually exclusive notions on a Cantor set, and follows from the Jewett-Krieger Theorem by 
\cite{Pe}.
\par
Third section concerns with uniform rigidity and measurable weak mixing. Author motivation for this topic is that
 a (nontrivial) measure-preserving weakly mixing transformation that is uniformly rigid would yield  an example of
 a measurable sensitive transformation that is not strongly measurably sensitive.
  For a subset $Y$ of a metric space $X$ and a measurable transformation  of $X$ authors of the paper \cite{JKLSS}
 define the notion of uniformly rigid transformation on $Y$ and prove Theorem 3.4 that is reminiscent of Egorov’s 
Theorem by P. Halmos \cite{Ha}.
In forth section authors present asymptotic convergence behavior. Let $X$ be a compact metric space and let $T$ be a
 finite measure-preserving ergodic transformation. Authors prove:
\begin{proposition} \label{T14}
  If $T$ is uniformly rigid, then
 the uniform rigidity sequence has zero density.
\end{proposition}
\par
The aim of section five is to study group action and generalized uniform rigidity.
   Let $G$ be a countable group endowed with the discrete topology acting faithfully on a finite measure space by
 measure-preserving transformations. Following authors of the paper \cite{JKLSS} the action of $G$ is uniformly
 rigid if there exists a sequence $\{g_{i}\}$ of group elements that leaves every compact $K \subset G,$ denoted
 $ g_{i} \to \infty,$ such that $d(x, g_{i} \cdot x) \to 0$ uniformly. The main result of the section is Theorem 5.3: 

\begin{theorem} \label{T15}
Let $X$ admit a weakly mixing group action and a uniformly rigid action by nontrivial subgroups of a fixed group of
 automorphisms $G.$ Then there exists a $G-$action on $X$ that is simultaneously weakly mixing and uniformly rigid.  
\end{theorem}

Authors formulate several interesting questions  that arise under
 investigations of weak mixing and uniform rigidity.

Some results and methods that are connected with topics of this and next section are considered in the paper \cite{Gl1}

\section{Actions of groups and semigroups}
\label{S5} \vspace{-4pt}
Furstenberg and Berent investigate the action of abelian semigroups on the
 torus ${\bf T}^d$ for $d=1$ and $d>1$ respectively. The authors of the paper \cite{BFLM} extend to the noncommutative case
 some results of Furstenberg and Berent. Author`s results answer problems raising by H. Furstenberg \cite{Fu1} and by Y. Guivarc{'}h
 [private communication to authors of the paper \cite{BFLM}].

Let $\nu$  be a probability measure on $SL_d ({\bf Z})$ satisfying the moment condition
 ${\bf E}_{\nu}(\parallel g \parallel^{\epsilon}) < \infty $ for some $\epsilon$. The authors of the paper \cite{BFLM} 
show that if the group generated by the support of  $\nu$ is large enough, in particular if this group is Zariski dense in 
$SL_d$, for any irrational $x \in {\bf T}^d$ the probability measures $\nu^{*n}*{\delta}_x$ tend to the uniform measure on 
${\bf T}^d.$  If in addition $x$ is Diophantine generic, authors show 
this convergence is exponentially fast.
\par
   This interesting survey includes resent results on rigidity theory by 
M. Einsiedler, E. Lindenstrauss \cite{EL} and
 by G.A. Margulis \cite{Ma1}, 
convolution of measures, on $\nu-$stiff action, on Fourier coefficients of measures and on notions of coarse dimension.
\par
Let the action of semigroup ${\Gamma }$ on ${\bf T}^d$ satisfy the following three conditions:  
($\Gamma\! -\! 0$) ${\Gamma } < SL_{d}({\bf R}),$ (${\Gamma }\! -\! 1$) ${\Gamma}$ acts strongly irreducibly 
on ${\bf R}^d,$  (${\Gamma}\! -\! 2$) ${\Gamma}$ contains a proximal element: there is some 
$g \in \Gamma$ with a dominant eigenvalue which is a simple root of its characteristic polynomial.
\par
   In Section 1 authors formulate main result of the paper. 
\begin{theorem} \label{T6}
 Let ${\Gamma} < SL_{d}({\bf R})$ satisfy 
($\Gamma\! -\! 1$) and  ($\Gamma\! -\! 2$) above, and let $\nu$ be a probability measure supported on a set of generators of 
${\Gamma }$ satisfying $\sum_{g \in \Gamma} \nu(g) \parallel g \parallel^{\epsilon} < \infty $ for some $\epsilon > 0.$
 Then for any $0 < \lambda < \lambda_{1}(\nu)$ there is a constant $C = C(\nu, \lambda)$ so that if for a point
 $x \in {\bf T}^d$ the measure $\mu_n = \nu^{*n}*{\delta}_x$ satisfies that for some 
$a \in  {{\bf Z}}^d \setminus \{0\} \; \mid \hat \mu_{n}(a) \mid > t > 0,$ with $n > C\cdot \ log (\frac{2\parallel a \parallel}{t}),$ 
then $x$ admits a rational approximation 
$p / q$ for $p \in {{\bf Z}}^d$ and $q \in {{\bf Z}}_{+}$ satisfying $\parallel x - \frac{p}{q} \parallel < \exp^{-\lambda n}$
 and    $\mid q  \mid <  (\frac{2\parallel a \parallel}{t})^C.$
\end{theorem}
Authors of \cite{BFLM} denote the theorem as Theorem A.
\par
Section 2 is devoted to the deduction of two corollaries from Theorem A. Let in the corollaries ${\Gamma}$ and $\nu$  be as 
in theorem A.
\begin{corollary} \label{C0}
 Let $x \in {{\bf T}^d} \setminus ({\bf Q}/{\bf Z})^d. $ Then the measures
 $\mu_n = \nu^{*n}*{\delta}_x$ converge to the Haar measure $m$ on ${\bf T}^d$ in weak-* topology. 
\end{corollary}
 This is authors \cite{BFLM} Corollary B.
Next corollery is the authors  \cite{BFLM} Corollary C:

\begin{corollary} \label{C5}
 Let $x \in {{\bf T}^d}$ and $\mu_n = \nu^{*n}*{\delta}_x.$ Then there are $c_1, c_2$ depending only on  
$\nu$ so that the following holds: (1) Assume $x$ is Diophantine generic in the sense that for some $M$ and
 $Q \; \parallel x - \frac{p}{q} \parallel > q^{-M}$ for all integers $q \ge Q$ and $p \in{\bf Z}^d.$ Then 
for $n > c_1 \log Q \;  \max_{b \in {\bf Z}^d , \parallel b \parallel < B}   \mid \hat \mu_{n}(b) \mid < B \exp^{{- c_2 n}/M}.$ 
(2) Assume $x \notin {\bf Q}/{\bf Z})^d.$ Then there is a sequence 
$n_{i} \to \infty$ along which  $\max_{b \in {\bf Z}^d , 0<\parallel b \parallel < \exp^{ c_2 n_{i}}}   \mid \hat \mu_{n}(b) \mid <  \exp^{- c_2 n_i}.$   
\end{corollary}

Section 3 gives the deduction of authors' solution  of Furstenberg problem from the authors \cite{BFLM}  Proposition 3.1: 
\begin{proposition} \label{T17}
Let ${\Gamma}$ and $\nu$  be as in theorem A, $0 < \lambda < \lambda_{1}(\nu).$ Then for some constant $C$ depending on
 $\nu , \lambda$ the following holds: for any probability measure $\mu_0$ on ${\bf Z}^d,$ if $\mu_n = \nu^{*n}*{\mu_0}$ 
has a nontrivial Fourier coefficient 
$a \in {\bf Z}^d \setminus \{0\} \; \mid \hat \mu_{n}(a) \mid > t,$ with $n > C \cdot \ log (\frac{2\parallel a \parallel}{t}),$ 
then $\mu_{0}(W_{Q, \exp^{-\lambda n}}) > (\frac{t}{2})^C$ where $Q = (\frac{2\parallel a \parallel}{t})^C.$
\end{proposition}
 Theorem A follows from Proposition 3.1.
\par
Section 4 is devoted to random matrix products. It includes estimates of the metric on ${\bf P}^{ d-1}$ and random walks.
In Section 5 two notions of coarse dimension are discussed.
  Section 6 describes the structure of the set of $t-$large Fourier coefficients.
The last Section "Granulated measures"  gives the prove of Proposition 3.1.
\par
The results of the paper\cite{BFLM}  will be of use to specialists interested in Diophantine approximation, 
measure theory and algebraic dynamics.

\section{Rigid cohomology}
\label{S6} \vspace{-4pt}

Let $p$ be a prime, $n$ a positive integer, and ${\bf F}_q$ the finite field with $q = p^n$ elements.
Let ${\bf Q}_q$ denote the unique unramified extension of degree $n$ of the field of $p$-adic numbers. 
Let $U$ be an open dense subscheme of the projective space ${\bf P}^{1}_{{\bf Q}_q}$ with nonempty complement $Z$.
     Let $V$ be the rigid analytic subspace of ${\bf P}^{1}_{{\bf Q}_q}$ which is the complement of the union of 
the open disks of radius $1$ around the points of $Z$. 
  A Frobenius structure on $ {\cal E}$ with respect to $\sigma$ is an isomorphism 
${\cal F}: \sigma^{*}{\cal E} \simeq {\cal E} $ of vector bundles with connection defined on some strict neighborhood of $V$.
\par
A meromorphic connection on ${\bf P}^{1}$ over a $p$-adic field admits a Frobenius structure defined over a suitable 
rigid analytic subspace. 
Authors of the paper\cite{KT}  give an effective convergence bound for this Frobenius structure by studying the effect of 
changing the Frobenius lift. They also give an example indicating that their bound is optimal.
\par
The techniques used are computational. This is a good place to see the interplay between matrix representation of a 
Frobenius structure and a Gauss-Manin connection.

The theory of rigid $p$-adic cohomology are developed by Berthelot \cite{Be} and others.
Rigid cohomology in some sense extends crystalline cohomology. Review of some novel results and applications of
 crystalline cohomology is given in paper \cite{Gl}.

\section{Superrigidity}
\label{S7} \vspace{-4pt} 

The notion of property (T) for locally compact groups was defined by 
D. Kazhdan \cite{Ka} and the notion of relative property (T)  
for inclusion of countable groups $\Gamma_0 \subset \Gamma$ was defined by 
G. Margulis  \cite{Ma2}.

The concept of superrigidity was introduced by G. D.  Mostow  \cite{Mo1}
and by G. A. Margulis  \cite{Ma3}
in the context of studying the structure of lattices in rank one and higher rank Lie groups respectively.
The first result on orbit equivalent (OE) superrigid actions  was obtained by A. Furman \cite{Fur},
 who combined the cocycle superrigidity by 
R. Zimmer \cite{Zi} with ideas from geometric group theory to show that the actions 
$SL_{n} ({\bf Z}) \to  T^n (n \ge 3)$ are OE superrigid. The deformable actions of rigid 
groups are OE superrigid by S. Popa \cite{Po}.

The paper \cite{Io} presents a new class of orbit equivalent superrigid actions.
The main result of the paper \cite{Io} is the Theorem A on orbit equivalence (OE) superrigidity. As a consequence of 
Theorem A the author can constructs uncountable many non-OE profinite actions for the arithmetic groups 
$SL_n ({\bf Z}) (n \ge 3)$, as well as for their finite subgroups,  and for the groups that are semi direct products of groups  
$SL_m ({\bf Z})$ and ${\bf Z}^m (m \ge 2)$. The author deduces Theorem A as a consequence of the 
Theorem B on cocycle superrigidity. 
\par
Let $\Gamma \to X$ be a free ergodic measure-preserving profinite action (i.e., an inverse 
limit of actions $\Gamma  \to X_n$ with $X_n$ finite) of a countable property (T) group 
$\Gamma $ (more generally, of a group $\Gamma$ which admits an infinite normal subgroup $\Gamma_0$ such 
that the inclusion $\Gamma_0 \subset \Gamma$ has relative property (T) and ${\Gamma} / {\Gamma_0}$  is finitely 
generated) on a standard probability space $X$. The author prove that if $\omega:\Gamma \times X \to \Lambda$ 
is a measurable cocycle with values in a countable group $\Lambda$, then $\omega$ is a cohomologous to a cocycle 
$\omega^{'}$ which factors through the map $\Gamma \times X \to \Gamma \times X_n$, for some $n$. As a corollary, 
he shows that any free ergodic measure-preserving action $\Lambda \to Y$ comes from a (virtual) 
conjugancy of actions.

\section{Conclusion}
\label{S8} \vspace{-4pt}

Rigid, hard and soft problems and results in arithmetic geometry have presented. 
Diverse notions of rigidity and respective novel results are reviewed.

\end{document}